\documentclass[12pt]{amsart}
\usepackage{amsmath}
\usepackage{amscd}
\usepackage{amssymb}
\usepackage{amsfonts}
\usepackage{amsthm}
\usepackage{bbm}
\usepackage{cancel}
\usepackage{color}
\usepackage{eucal}
\usepackage{enumerate,yfonts}
\usepackage{enumitem}

\usepackage{pdfsync}
\usepackage[all,cmtip]{xy}
\usepackage{graphicx}
\usepackage{graphics}
\usepackage{hyperref}
\usepackage{latexsym}
\usepackage{mathrsfs}
\usepackage{placeins}
\usepackage{pstricks}
\usepackage{bookmark}
\usepackage{mathtools}
\usepackage{stmaryrd}
\usepackage{url}

\newtheorem{thm}{Theorem}[section]

\newtheorem{lemma}[thm]{Lemma}

\newtheorem{thm-dfn}[thm]{Theorem-Definition}

\newtheorem{remark}[thm]{Remark}

\newtheorem{conv}[thm]{Convention}

\numberwithin{equation}{section}

\newcommand{\nc}{\newcommand}

\newcommand{\cO}{\mathcal{O}}

\newcommand{\cF}{\mathcal{F}}

\newcommand{\cH}{\mathcal{H}}

\newcommand{\cP}{\mathcal{P}}
\newcommand{\cL}{\mathcal{L}}

\newcommand{\cX}{\mathcal{X}}

\newcommand{\bC}{{\mathbb C}}

\newcommand{\bR}{{\mathbb R}}
\newcommand{\bH}{{\mathbb H}}

\newcommand{\bi}{{\mathbf i}}

\newcommand{\Lt}{{\mathfrak{t}}}
\newcommand{\Lg}{{\mathfrak g}}

\newcommand{\rB}{\mathrm{B}}

\newcommand{\binv}{{\backslash\hspace{-0.04in}\backslash}}

\nc{\on}{\operatorname}
\nc{\codim}{{\operatorname{codim}}}
\nc{\img}{{\operatorname{Im}}}
\nc{\IC}{{\operatorname{IC}}}

\nc{\lp}{{\left(}}
\nc{\rp}{{\right)}}

\newcommand{\beqn}{\begin{equation*}}
\newcommand{\eeqn}{\end{equation*}}

\newcommand{\beq}{\begin{equation}}
\newcommand{\eeq}{\end{equation}}

\newcommand{\bern}{\begin{eqnarray*}}
\newcommand{\eern}{\end{eqnarray*}}

\newcommand{\WquotMap}{g}

\begin{document}

\title[Errata and Notes]{Errata and Notes on the Paper ``A Generalization of
                                     Springer Theory Using Nearby Cycles"}

\author{Mikhail Grinberg}

\date{February 26, 2020}

\maketitle

\begin{abstract}
We provide some corrections and clarifications to the paper \cite{Gr3} of the title.
In particular, we clarify the ``left/right'' conventions on complex reflection groups
and their braid groups.  Most importantly, we fill in a gap related to the treatment
of cuts in the Picard-Lefschetz theory part of the argument.  The statements of
the main results are not affected.
\end{abstract}

\tableofcontents

\section{Introduction}

The paper \cite{Gr3} of the title contains a gap in the proof of its main 
result, \cite[Theorem 3.1]{Gr3}, related to the treatment of cuts in the
Picard-Lefschetz theory part of the argument.  Namely, the proofs of
\cite[Lemmas 4.2 $\&$ 4.3]{Gr3} are not satisfactory as written, and to the
author's knowledge, can not be fixed without substantial further argument.
The main goal of this document is to fill in this gap.  This is done by slighltly
modifying the statements of \cite[Lemmas 4.2 $\&$ 4.3]{Gr3}, and by providing
proofs of the modified lemmas.  In addition, we take this opportunity to provide
some minor corrections, clarifications, and additional details for the rest of
\cite{Gr3}.

In more detail, the contents of this document are as follows.  In Section
\ref{sec-left-right}, we clarify our ``left/right'' conventions on complex reflection
groups and their braid groups.  The material here is mostly notational, and
some readers may find the notation to be overkill.  However, the author has 
found it easy to make mistakes or to create ambiguities related to such
``left/right'' conventions (starting with the convention for multiplying loops in
a fundamental group), and hopes that the notation presented here is helpful
in avoiding such mistakes and ambiguities.  In Section
\ref{sec-cyclic-existence}, we describe in detail and fix the gap related
to \cite[Lemmas 4.2 $\&$ 4.3]{Gr3}.  The problem is described in Section
\ref{subsec-problem}.  Statements of the modified lemmas are given in
Section \ref{subsec-modifying}, and the proofs are given in Sections
\ref{subsec-proof-4p2}-\ref{subsec-proof-4p3}.  Sections
\ref{subsec-stable-case}-\ref{subsec-nonstable-case} indicate how
to adapt the rest of the proof of \cite[Theorem 3.1]{Gr3} to the modified
versions of the lemmas.  In addition, Section \ref{subsec-stable-case}
corrects a couple of unrelated minor issues, and Section
\ref{subsec-nonstable-case} substantially expands on the corresponding
\cite[Section 4.3]{Gr3} by providing significant further detail.  Finally,
in Section \ref{sec-sec-5-6}, we make several corrections to
\cite[Sections 5-6]{Gr3}, and remark on the relation of the material in
\cite[Section 6]{Gr3} to the recent paper \cite{GVX}.

We number the bibliography to extend the numbering in \cite{Gr3}.  We also 
provide the publication data for the references \cite{Gr1} and \cite{Gr2}.

\section{Left/Right Conventions}
\label{sec-left-right}

We begin by clarifying some basic conventions, pertaining to complex reflection
groups and their braid groups, that are not made sufficiently clear in \cite{Gr3}.

\subsection{The Centralizer of a Complex Reflection Group}
\label{subsec-centralizer}

Let $V$ be a finite dimensional complex vector space, and let $W \subset GL(V)$
be a finite complex reflection group acting on $V$.  See \cite{Bro} for an introduction
to such groups, also known as Shephard-Todd groups.  We say that $v \in V$ is
regular if the stabilizer of $v$ in $W$ is trivial.  We denote by $V^{reg} \subset V$
the set of all regular $v \in V$.  Pick a $W$-orbit $\cO \subset V^{reg}$.  Let
$S(\cO)$ be the group of permutations of the finite set $\cO$.  Note that $W$ is
naturally a subgroup of $S(\cO)$.  We denote by $W^\sharp (\cO) \subset S(\cO)$
the centralizer of $W$ in $S(\cO)$.

The group $W^\sharp (\cO)$ is isomorphic to $W$ as an abstract group; but we
prefer to think of $W^\sharp (\cO)$ as isomorphic to the opposite group $W^{op}$
of $W$.  More precisely, for every $v \in \cO$, we can define a bijection of sets:
\beqn
\varphi_v : W^\sharp (\cO) \to W,
\eeqn
by the following condition:
\beqn
\varphi_v (a) (v) = a (v) \;\; \text{for every} \;\; a \in W^\sharp (\cO).
\eeqn
It is not hard to check that $\varphi_v : W^\sharp (\cO) \to W$ is an
anti-homomorphism of groups.  Indeed, for every pair $a, b \in W^\sharp (\cO)$,
we have:
\beqn
\varphi_v (a b) (v) = (a b) (v) 
                             = a (\varphi_v (b) (v)) 
                             = \varphi_v(b) (a (v))
                             = (\varphi_v(b) \varphi_v(a)) (v),
\eeqn
where the third equality holds because $a \in W^\sharp (\cO)$ and
$\varphi_v(b) \in W$ commute as elements of $S(\cO)$, by the definition of 
$W^\sharp (\cO)$.  Thus, we obtain a group isomorphism: 
\beqn
\varphi_v : W^\sharp (\cO) \to W^{op}.
\eeqn

Let $v_1, v_2 \in \cO$ be a pair of points, and let $a \in W^\sharp (\cO)$
be the unique element such that $v_2 = a(v_1)$.  One can check that the
composition:
\beqn
\varphi_{v_2}^{-1} \circ \varphi_{v_1} : W^\sharp (\cO) \to W^\sharp (\cO),
\eeqn
is given by:
\beqn
\varphi_{v_2}^{-1} \circ \varphi_{v_1} : b \mapsto a b a^{-1}.
\eeqn
Thus, we see that the group $W^\sharp (\cO)$ is non-canonically isomorphic
to $W^{op}$.

\begin{remark}  
The author is not aware of any source in the literature that highlights the group
$W^\sharp (\cO)$, as distinct from $W$ or $W^{op}$.  Of course, any statement
about $W^\sharp (\cO)$ can be expressed in terms of $W$ and $W^{op}$.
However, we feel that some things are made clearer by highlighting the group
$W^\sharp (\cO)$.  

Note that the action of $W^\sharp (\cO)$ on $\cO$ does not, in general, extend to
a linear action on $V$.  In fact, it does not even extend to a continuous action on
$V^{reg}$.  If we pick a point $v \in \cO$, to make the identifications $\cO \cong W$,
$W^\sharp (\cO) \cong W^{op}$, the action of $W^\sharp (\cO)$ on $\cO$ becomes
the right action of $W$ on itself.  However, this point of view obscures the geometric
distinction between the non-linear action of $W^\sharp (\cO)$ and the linear action
of $W$ on $\cO$.
\end{remark}

\subsection{The Braid Group of a Complex Reflection Group}
\label{subsec-braid-group}

Let $Q = W \binv V$.  Recall that, as an algebraic variety, the quotient $Q$ is
isomorphic to $V$ (see \cite[Theorem 4.1]{Bro}).  We will not be using any explicit
linearization of $Q$; but we note that it naturally possesses an origin, given by the
orbit of $0 \in V$.  Let $\WquotMap : V \to Q$ be the quotient map, and let
$Q^{reg} = \WquotMap (V^{reg}) \subset Q$.  Pick a point $q \in Q^{reg}$.
The braid group $B_W (q)$ of $W$ at $q$ is defined as a fundamental group:
\beqn
B_W (q) = \pi_1 (Q^{reg}, q).
\eeqn
At this point, it is essential to state the following convention, which ensures that
monodromy representations of fundamental groups are, in fact, left representations
(see \cite[Remarks 2.3.1, 2.6.3]{Sz}).

\begin{conv}
We multiply loops in a fundamental group $\pi_1 (X, x)$ by tracing the second
loop first.
\end{conv}

Note that $\WquotMap$ restricts to a covering map:
\beqn
\WquotMap^{reg} : V^{reg} \to Q^{reg},
\eeqn
and that $W$ acts on $V^{reg}$ by deck transformations.  Therefore, writing 
$\cO_q = \WquotMap^{-1} (q)$, we obtain a group homomorphism:
\beqn
\eta_q : B_W (q) \to W^\sharp (\cO_q),
\eeqn
given by the monodromy of the covering map $\WquotMap^{reg}$.  It is not hard to
check that the homomorphism $\eta_q$ is surjective.  Thus, we can think of
$W^\sharp (\cO_q)$ as a quotient of $B_W (q)$.  If we further pick a point
$v \in \cO_q$, we can define a group homomorphism:
\beqn
\eta^\flat_v = I_W \circ \varphi_v \circ \eta_q : B_W (q) \to W,
\eeqn
where $I_W : W^{op} \to W$ is the inverse map (we make no distinction
between the group homomorphism $I_W$ and its inverse, as the underlying
sets of $W$ and $W^{op}$ are the same).

\begin{remark}\label{locsysoverQ}
The groups $W^\sharp (\cO_q)$, for all $q \in Q^{reg}$, naturally form a local
system over $Q^{reg}$.  The same is true of the groups $B_W (q)$.  The maps
$\eta_q$, for all $q \in Q^{reg}$, give a map between these local systems.
\end{remark}

\subsection{The Case of a Coxeter Group}
\label{subsec-coxeter-case}

Suppose now that $V_\bR$ is a finite dimensional real vector space, that $W$ is
a finite Coxeter group acting on $V_\bR$, and that $V$ is the complexification of
$V_\bR$.  All of the above constructions can be applied to the pair $(V, W)$.
We let $V^{reg}_\bR = V^{reg} \cap V_\bR$, $Q_\bR = \WquotMap (V_\bR)$, and
$Q^{reg}_\bR = \WquotMap (V^{reg}_\bR)$.

\begin{remark}\label{wsharpR}
The space $Q^{reg}_\bR$ is contractible.  Therefore, by Remark \ref{locsysoverQ},
the groups $W^\sharp (\cO_q)$, for all $q \in Q^{reg}_\bR$, are naturally isomorphic
to each other.  The same is true of the groups $B_W (q)$, for all $q \in Q^{reg}_\bR$.
Based on this, we can unambiguously write $W^\sharp = W^\sharp (\cO_q)$,
$B_W = B_W (q)$, and $\eta = \eta_q : B_W \to W^\sharp$, where $q$ is
understood to be some point of $Q^{reg}_\bR$.
\end{remark}

\begin{remark}\label{varphiR}
Let $V^+ \subset V_\bR$ be a connected component of $V^{reg}_\bR$.  Then,
for every $v_1, v_2 \in V^+$, we have:
\beqn
\varphi_{v_1} = \varphi_{v_2} : W^\sharp \to W.
\eeqn
Based on this, we can unambiguously write:
\beqn
\varphi^{}_{V^+} = \varphi_v : W^\sharp \to W
\;\; \text{and} \;\;
\eta^\flat_{V^+} = \eta^\flat_v : B_W \to W,
\eeqn
where $v$ is understood to be some point of $V^+$.
\end{remark}

Of special interest is the case of the Weyl group of a complex semisimple Lie algebra
$\Lg$.  Every Cartan subalgebra $\Lt \subset \Lg$ comes equipped with a canonical
real form $\Lt_\bR \subset \Lt$, such that $\bi \cdot \Lt_\bR$ generates a compact torus
(where $\bi = \sqrt {-1}$).  Associated to each $\Lt \subset \Lg$, we have a Weyl group
$W [\Lt]$, which is a finite Coxeter  group acting on $\Lt_\bR$.  By Remark
\ref{wsharpR}, we also have the groups $W^\sharp [\Lt]$ and $B_{W [\Lt]}$,
associated to $W [\Lt]$.

\begin{remark}\label{canonicalWsharp}
The groups $W^\sharp [\Lt]$, for different choices of $\Lt$, are naturally isomorphic 
to each other.  The same is true of the braid groups $B_{W [\Lt]}$, but not of the Weyl
groups $W [\Lt]$.  Based on this, we can unambiguously write $W^\sharp[\Lg] =
W^\sharp[\Lt]$ and $B_{W [\Lg]} = B_{W [\Lt]}$, where $\Lt$ is understood to be some
Cartan subalgebra of $\Lg$.  Furthermore, the homomorphism $\eta : B_{W [\Lg]} \to
W^\sharp[\Lg]$ of Remark \ref{wsharpR} is independent of the Cartan subalgebra
used to define it.
\end{remark}

\subsection{Comments on \cite[Section 2.1]{Gr3}}
\label{subsec-comments-2.1}

In \cite[Section 2.1]{Gr3}, we gave a brief summary of the approach to constructing the
Springer representations due to Lusztig and Borho-MacPherson.  Here, we connect
the exposition in that section to the terminology and notation introduced in Sections
\ref{subsec-centralizer}-\ref{subsec-coxeter-case} above.  We hope that the use of the
centralizer group $W^\sharp = W^\sharp [\Lg]$ in this context will add some clarity, as
well as facilitate the connection to the nearby cycles approach to Springer theory (see
Section \ref{subsec-comments-2.2.1} below).

At the top of \cite[p. 414]{Gr3}, we consider the finite covering map:
\beqn
q^{rs} : \tilde \Lg^{rs} \to \Lg^{rs}.
\eeqn
After making a choice of a point $x \in \Lg^{rs}$ and a positive Weyl chamber for the
Cartan subalgebra $\Lg_x \subset \Lg$ containing $x$, we claim that the Weyl group
$W = W[\Lg_x]$ of $\Lg$ acts as the deck transformations of the covering $q^{rs}$.

The same idea can be expressed somewhat more precisely, by first making the
following observation.

\begin{remark}\label{wsharpAsDeck}
The group of deck transformations of the covering map $q^{rs}$ is canonically
isomorphic to the group $W^\sharp = W^\sharp [\Lg]$ of Remark \ref{canonicalWsharp}.
This can be seen by identifying the fiber $q^{-1} (x)$ with the set of all Weyl chambers
in $\Lg_x$, i.e., of all connected components of $(\Lg_x)^{reg}_\bR$.
\end{remark}

Let $\Lg^+_x \subset \Lg_x$ be the chosen positive Weyl chamber.  In order to
obtain an action of the Weyl group $W = W [\Lg_x]$ as the deck transformations
of $q^{rs}$, we must use Remark \ref{wsharpAsDeck} plus the homomorphism:
\beqn
\varphi_{\Lg^+_x}^{-1} \circ I_W : W \to W^\sharp,
\eeqn
where $I_W : W \to W^{op}$ is the inverse map, and $\varphi_{\Lg^+_x} :
W^\sharp \to W^{op}$ is the homomorphism of Remark \ref{varphiR}.  

Thus, one can say that Lusztig's construction, described in the paragraph following
\cite[Proposition 2.1]{Gr3}, naturally produces an action of the centralizer group
$W^\sharp$ on the push-forward sheaf $P$; and that we need to make some auxiliary
choices, and to apply the inverse map, in order to convert this action into an action of
the Weyl group $W$.  In particular, we can restate \cite[Theorem 2.2 (ii)]{Gr3} as
follows.

\begin{remark}\label{wsharponP}
The action of $W^\sharp$ on $P$, arising from Remark \ref{wsharpAsDeck}
via Lusztig's construction, gives an isomorphism $\bC [W^\sharp] \cong
\on{End} (P)$.
\end{remark}

\subsection{Comments on \cite[Section 2.2.1]{Gr3}}
\label{subsec-comments-2.2.1}

At the top of \cite[p. 416]{Gr3}, we consider the adjoint quotient map
$f : \Lg \to W \binv \Lt$.  Here, we clarify the definition of the nearby cycles
sheaf $P_f$ of $f$ and the statement of \cite[Theorem 2.5 (ii)]{Gr3}.  Let us
write $Q = W \binv \Lt$ for the target of $f$.  In order to fully specify the
sheaf $P_f$, we must fix a test-arc $\gamma$ in $Q$.  By this we mean
an embedded complex analytic arc $\gamma : U \to Q$, where $U$ is a 
neighborhood of zero in $\bC$, such that $\gamma (0) = 0$, and
$\gamma (z) \in Q^{reg}$ for $z \neq 0$; see the discussion at
\cite[p. 415]{Gr3}.  We will say that a test-arc $\gamma : U \to Q$ is real
if $\gamma (z) \in Q_\bR (= f(\Lt_\bR))$ for every $z \in U \cap \bR$.

\begin{remark}\label{realArcs}
Let $\rB_\epsilon \subset Q$ be a small ball around the origin, defined using
some complex analytic local coordinates on $Q$.  The local fundamental
group $\pi_1 (Q^{reg}_\bR \cap \rB_\epsilon)$ is trivial.  Therefore, the nearby
cycles sheaves $P_{f_\gamma}$, for all real test-arcs $\gamma$ in $Q$, are
naturally isomorphic to each other.
\end{remark}

Using Remark \ref{realArcs}, we can resolve the ``up-to-isomorphism'' ambiguity
in the definition of $P_f$ by assuming that $P_f = P_{f_\gamma}$ for some real
test-arc $\gamma$ in $Q$.  With this convention, the monodromy group acting
on $P_f$ is the braid group $B_{W [\Lg]}$ of Remark \ref{canonicalWsharp}.
To keep the notation close to \cite[Section 2.2.1]{Gr3}, we will write $B_W =
B_{W [\Lg]}$.  We can now restate \cite[Theorem 2.5 (ii)]{Gr3} as follows.

\begin{remark}\label{bwonP}
The monodromy action $\mu : B_W \to \on{Aut} (P_f)$ factors through the
homomorphism $\eta : B_W \to W^\sharp$ of Remark \ref{canonicalWsharp},
producing an action of $W^\sharp$ on $P_f$.  The isomorphism of
\cite[Theorem 2.5 (i)]{Gr3} agrees with the actions of $W^\sharp$ on both
sides (see Remark \ref{wsharponP} for the action of $W^\sharp$ on the LHS).
\end{remark}

Remark \ref{bwonP} seems like the best way to clarify the meaning of
\cite[Theorem 2.5 (ii)]{Gr3}.  However, in order to read this theorem as
stated (i.e., in terms of the Weyl group $W$ rather than the centralizer
group $W^\sharp$), we must clarify what is meant by the ``natural
homomorphism $B_W \to W$''.  For this, we need to invoke the same
auxiliary choices that were used in Section 2.4 above.  Namely, a Cartan
subalgebra $\Lg_x \subset \Lg$, and a positive Weyl chamber
$\Lg^+_x \subset \Lg_x$.  With these choices, we assume that
$W = W[\Lg_x]$, and we use the homomorphism:
\beqn
\eta^\flat_{\Lg^+_x} : B_W \to W,
\eeqn
of Remark \ref{varphiR}, as the ``natural homomorphism $B_W \to W$''
of \cite[Theorem 2.5 (ii)]{Gr3}.

Thus, we see that both from the point of view of Lusztig's construction, and
from the point of view of the nearby cycles construction, the centralizer group
$W^\sharp = W^\sharp [\Lg]$ emerges as the natural group of symmetries
of the perverse sheaf $P$.

\section{Existence of a Cyclic Picard-Lefschetz Class}
\label{sec-cyclic-existence}

Our comments in Sections \ref{subsec-comments-2.1} and
\ref{subsec-comments-2.2.1} above amounted to clarifications rather than
corrections.  In this section, we address an actual gap in the proof of
\cite[Theorem 3.1]{Gr3}.  We being by explaining the problem.

\subsection{The Problem}
\label{subsec-problem}

The proof of \cite[Theorem 3.1]{Gr3} contains a gap in the part of the
argument dealing with Picard-Lefschetz cuts.  More precisely, the proofs of
\cite[Lemmas 4.2 $\&$ 4.3]{Gr3} are not satisfactory as written.  To explain
the problem, we focus on the proof of \cite[Lemma 4.2]{Gr3}.

Before proceeding to the substantive issue, we need to clarify the first sentence
of this proof: ``For each $w \in W$, pick a lift $b_w \in B_W$.''  What is meant
here is that, in the notation of Sections
\ref{subsec-centralizer}-\ref{subsec-braid-group}
above, we have a surjective group homomorphism:
\beqn
\varphi^{}_{e_0} \circ \eta^{}_\lambda : B_W \to W^{op},
\eeqn
and the lift $b_w$ is taken with respect to this homomorphism. 

With this understanding, in the proof of \cite[Lemma 4.2]{Gr3}, we proceed
by constructing a family of Picard-Lefschetz classes: 
\beqn
\{ \mu (b_w) \, u_0 = u (e_w, \gamma_w, \cO_w) \}_{w \in W} \, ,
\eeqn
parametrized by the Weyl group $W$ of the polar representation $G | V$,
associated to the Cartan subspace $c$.  Here, each $e_w$ is a critical
point of the algebraic function:
\beqn
l |_F : F \to \bC,
\eeqn
each $\gamma_w$ is a smooth path in $\bC$, connecting $l (e_w)$ to a fixed,
large $\xi_0 > 0$, and each $\cO_w$ is an orientation of the ``positive'' real
subspace $T_{e_w} [\gamma_w] \subset T_{e_w} F$, corresponding to the
path $\gamma_w$.  

Next, we observe that the points $\{ e_w \}_{w \in W}$ are the only critical
points of the function:
\beqn
\hat l |_{\hat F} : {\hat F} \to \bC,
\eeqn
where $\hat F$ is a certain compactification of $F$ relative to $l$, and
$\hat l |_{\hat F}$ is the extension of $l |_F$ to $\hat F$.  We then conclude
that the classes $\{ \mu (b_w) \, u_0 \}_{w \in W}$ form a basis of the relative
homology group:
\beqn
H^{-d}_l (\cF \, P) \cong H_{d-r} (F, \{ \xi(y) \geq \xi_0 \}; \bC),
\eeqn
where $\xi = \on{Re} (l) : F \to \bR$.  It is this conclusion that does not
seem justified, because of insufficient control over the paths
$\{ \gamma_w \}_{w \in W}$.  More precisely, each of the $\gamma_w$
is guaranteed to satisfy conditions (i)-(iv) of \cite[p. 424]{Gr3}.  If, in addition,
we knew that:

(v) $\gamma_{w_1} (t_1) \neq \gamma_{w_2} (t_2)$, for all $w_1 \neq w_2$
and all $(t_1, t_2) \in [0, 1] \times [0, 1] \setminus \{ (1, 1) \}$,

\noindent
then the conclusion that $\{ \mu (b_w) \, u_0 \}_{w \in W}$ is a basis would 
be justified and standard.  However, our construction does not guarantee
condition (v).  It is easy to reduce \cite[Lemma 4.2]{Gr3} to the case where
the critical values $\{ l (e_w) \}_{w \in W}$ are distinct.  But nothing in our
construction ensures that the interiors of the paths
$\{ \gamma (e_w) \}_{w \in W}$ will not intersect.  This is a puzzling
oversight; but the author is unable to see how to complete the proof of
\cite[Lemma 4.2]{Gr3} without some substantive further argument.  In fact,
it is not clear to the author whether the lemma is true as stated.  The proof
of \cite[Lemma 4.3]{Gr3} is analogous, and suffers from the same problem.

\subsection{Modifying the Statements of \cite[Lemmas 4.2 $\&$ 4.3]{Gr3}}
\label{subsec-modifying}

The statements of \cite[Lemmas 4.2 $\&$ 4.3]{Gr3} refer to a particular
Picard-Lefschetz class $u_0 \in H_l^{-d} (\cF \, P)$.  However, in defining
this class (see \cite[p. 424]{Gr3}), we placed no restrictions on the critical
point $e_0 \in Z$, and no restrictions the path $\gamma_0 : [0,1] \to \bC$,
connecting $l(e_0)$ to $\xi_0$, in addition to the general conditions (i)-(iv)
of \cite[p. 424]{Gr3}.  On the other hand, the proof of \cite[Theorem 3.1]{Gr3}
made only some very mild use of the flexibility to make specific choices of
$e_0$ and $\gamma_0$.  For this reason, in order to fix the problem
described in Section 3.1 above, it seems easiest to put the following two
additional restrictions on the pair $(e_0, \gamma_0)$.  Recall that we write
$\xi = \on{Re} (l)$.

(A1)  For every $e \in Z$, we have $\xi (e) \leq \xi (e_0)$.

(A2) For every $t \in (0, 1]$, we have $\xi (\gamma_0 (t)) > \xi (e_0)$.

\begin{lemma}\label{modified4.2}
The statement of \cite[Lemma 4.2]{Gr3} holds with the additional assumption
that the pair $(e_0, \gamma_0)$ satisfies conditions {\em (A1)} and {\em (A2)}
above.
\end{lemma}

\begin{lemma}\label{modified4.3}
The statement of \cite[Lemma 4.3]{Gr3} holds with the additional assumption
that the pair $(e_0, \gamma_0)$ satisfies conditions {\em (A1)} and {\em (A2)}
above.
\end{lemma}

We will prove Lemmas \ref{modified4.2} and \ref{modified4.3} in Sections 3.3
and 3.4, respectively.  In Section 3.5, will indicate how to modify the proof of
\cite[Theorem 3.1]{Gr3} to make use of Lemmas \ref{modified4.2} and
\ref{modified4.3} in place of \cite[Lemmas 4.2 $\&$ 4.3]{Gr3}.

\subsection{Proof of Modified \cite[Lemma 4.2]{Gr3}}
\label{subsec-proof-4p2}

In this subsection, we provide a proof of Lemma \ref{modified4.2}.

{\em Step 1.}  Recall the dual Cartan subspace $c^* \subset V^*$ of $c \subset V$
(see \cite[Proposition 2.13]{Gr3}).  There exists a covector $l_1 \in c^*$, such that
$\on{Re} (l_1 (e)) < \on{Re} (l_1 (e_0))$ for every $e \in Z \setminus \{ e_0 \}$.
This follows form the fact that the finite set $Z$ is the set of vertices of its convex
hull.

{\em Step 2.}  By considering the straight line path from $l$ to $l+l_1$, and
replacing $l$ by $l+l_1$ in the statement of the lemma, if necessary, we can
assume that the following stronger version of condition (A1) holds:

(A3)  For every $e \in Z \setminus \{ e_0 \}$, we have $\xi (e) < \xi (e_0)$.

{\em Step 3.}  By further perturbing the covector $l \in c^*$, we can ensure that
the following two additional conditions hold.  We write $\zeta = \on{Im} (l)$.

(A4)  For every pair of distinct points $e_1, e_2 \in Z$, we have
$\xi (e_1) \neq \xi (e_2)$.

(A5)  For every pair of distinct points $e_1, e_2 \in Z$, we have
$\zeta (e_1) \neq \zeta (e_2)$.

{\em Step 4.}  Let $n$ be the order of $W$.  We index the elements of $W$ by
the set $I = \{ 0, \dots, n-1 \}$, so that the sequence $\{ \xi (w_i \, e_0) \}_{i \in I}$
is strictly decreasing.  This is possible by condition (A4) above.  Note that, by
condition (A3), we have $w_0 = 1 \in W$.  To avoid cumbersome notation, we
will write $e_i =  e_{w_i}$ for every $i \in I$.

{\em Step 5.}  In the proof of \cite[Theorem 3.1]{Gr3}, we used Picard-Lefschetz
classes defined using paths in $\bC$, connecting critical values of $l$ to the
fixed large $\xi_0 > 0$.  For the present argument, as a matter of convenience,
we prefer to use paths connecting critical values of $l$ to the vertical line 
$\{ z \in \bC \; | \; \on{Re} (z) = \xi_0 \}$.  For every $x \in \bC$, define a
``horizontal'' path $\gamma^{}_h [x] : [0, 1] \to \bC$ by:
\beqn
\gamma^{}_h [x] : t \mapsto (1-t) \cdot \on{Re} (x) + t \cdot \xi_0 +
\bi \cdot \on{Im} (x).
\eeqn
All Picard-Lefschetz classes in the argument that follows will be defined using
such horizontal paths.  Note that the paths $\{ \gamma^{}_h [l(e_i)] \}_{i \in I}$
avoid each other by condition (A5) above.

{\em Step 6.}  For each $i \in I$, consider the path:
\beqn
\gamma^{}_h [i] = \gamma^{}_h [l(e_i)] : [0, 1] \to \bC.
\eeqn
Recall the Hessian $\cH_{e_i} : T_{e_i} F \to \bC$ of $l |_F $ at $e_i$, and consider
the positive eigenspace $T_h [i] \subset T_{e_i} F$ of the non-degenerate real
quadratic form $\on{Re} (\cH_{e_i}) : T_{e_i} F \to \bR$ (here we use the Hermitian
inner product $\langle \; , \, \rangle$ of \cite[Proposition 2.12]{Gr3}).  Pick an
orientation $\cO_h [i]$ of $T_h [i]$.  The triple $(e_i, \gamma^{}_h [i], \cO_h [i])$
defines a Picard-Lfschets class $u [i] \in H^{-d}_l (\cF \, P)$, as in \cite[p. 424]{Gr3}.

{\em Step 7.}  By assumptions (A1) and (A2) of the lemma, we have
$u_0 = \pm u [0]$.  By reversing the orientation $\cO [0]$, if necessary, we can
assume that the sign is a plus: $u_0 = u [0]$.

{\em Step 8.}  Since $\{ e_i \}_{i \in I}$ are the only critical points of the function
$\hat l |_{\hat F}$, and since the paths $\{ \gamma^{}_h [i] \}_{i \in I}$ avoid each
other (see Step 5), the elements $\{ u [i] \}_{i \in I}$ form a basis of
$H^{-d}_l (\cF \, P)$.  This proves the assertion of the lemma about the
dimension of $H^{-d}_l (\cF \, P)$.

{\em Step 9.}  Let $M^0 \subset H^{-d}_l (\cF \, P)$ be the linear subspace
generated by the image of $u [0]$ under the monodromy action:
\beqn
\mu_* : B_W \to \on{End} (H^{-d}_l (\cF \, P)).
\eeqn
To complete the proof of the lemma, it suffices to show that $u [i] \in M^0$ for every
$i \in I$.  We will do so by induction on $i$.  For $i = 0$, the statement is obvious.
Let $I_+ =  I \setminus \{ 0 \}$.  Fix a $k \in I_+$, and assume that $u [i] \in M^0$ for
every $i \in I$ with $i < k$.  We need to show that $u [k] \in M^0$.

{\em Step 10.}  Recall the regular part $c^{reg} = c \cap V^{rs}$ of $c$.  We have
the following claim.

\noindent
{\bf Claim:} {\em
There exists a smooth path:
\beqn
\beta : [0, 1] \to c^{reg},
\eeqn
with $\beta (0) = e_0$ and $\beta (1) = e_k$, satisfying conditions {\em (B1)} and
{\em (B2)} below.  For each $i \in I$, let $\beta [i] : [0, 1] \to c^{reg}$ be the translate 
of $\beta$ by the element $w_i \in W$.  More precisely, we let:
\beqn
\beta [i] : t \mapsto w_i \, \beta [0] (t). 
\eeqn

{\em (B1)}  For each $i \in I$, the derivative $(\xi \circ \beta [i])'$ does not vanish at
any point of $[0, 1]$.

{\em (B2)}  For each pair $i, j \in I$, there is at most one $t \in [0, 1]$, such that:
\beqn
(\zeta \circ \beta [i]) (t) = (\zeta \circ \beta [j]) (t).
\eeqn
Moreover, the root $t$ satisfies $t \in (0, 1)$ and:
\beqn
(\zeta \circ \beta [i])' (t) \neq (\zeta \circ \beta [j])' (t).
\eeqn
}

\noindent
{\bf Proof:}  The path $\beta$ can be obtained as a $C^1$-small perturbation of the
straight line the path $\beta_{st} : [0, 1] \to c$, given by:
\beqn
\beta_{st} : t \mapsto (1-t) \cdot e_0 + t \cdot e_k \, .
\eeqn
Conditions (B1) and (B2) for $\beta_{st}$ in place of $\beta$ follow from conditions
(A4) and (A5) of Step 3.  It remains to note that (B1) and (B2) are $C^1$-open
conditions, and that the complement $c \setminus c^{reg}$ is a complex hyperplane
arrangement in $c$ (see \cite[Theorem 4.7]{Bro}).
\qed

{\em Step 11.}  Recall the quotient map $f : V \to Q$.  Consider the composition:
\beqn
\bar \beta = f \circ \beta : [0, 1] \to Q^{reg}.
\eeqn
Note that $\bar \beta (0) = f (e_0) = \lambda$ and
$\bar \beta (1) = f (e_k) = \lambda$.  Thus, $\bar \beta$ represents an element of
$B_W = \pi (Q^{reg}, \lambda)$.  We will show that $u [k] \subset M^0$ 
(see Step 9 above) by examining the monodromy operator:
\beqn
\mu_* (\bar \beta) : H^{-d}_l (\cF \, P) \to H^{-d}_l (\cF \, P).
\eeqn

{\em Step 12.}  To analyze the monodromy operator $\mu_* (\bar \beta)$, we
need to consider a parametrized version of the basis $\{ u [i] \}_{i \in I}$ of
Step 8.  We begin by defining a parametrized version of the stalk cohomology
group $H^{-d}_l (\cF \, P)$.  Let us write $T = [0, 1]$ (not to be confused with
the tangent spaces that appear in Step 6).  For every $t \in T$, let 
$F_t = f^{-1} (\bar \beta (t))$, and let $P_t \in \cP_G (E)$ be the nearby cycles
sheaf (with constant coefficients) given by the specialization of $F_t$ to
$\on{As} (F_t) = E$ (see the discussion at the top of \cite[p. 421]{Gr3}).  To
streamline the notation, we write:
\beqn
M_t = H^{-d}_l (\cF \, P_t).
\eeqn
Note that the vector spaces $\{ M_t \}_{t \in T}$ form a local system over
the segment $T$, and that we have $M_0 = M_1 = H^{-d}_l (\cF \, P)$.
For each pair $t_1, t_2 \in T$, we write:
\beqn
\nu_{t_1, t_2} : M_{t_1} \to M_{t_2} \, ,
\eeqn
for the parallel transport of the local system $\{ M_t \}_{t \in T}$.  Note that,
with this notation, we have:
\beqn
\nu_{0, 1} = \mu_* (\bar \beta) : M_0 \to M_1 \, .
\eeqn

{\em Step 13.}  
Next, we consider parametrized versions of the constructions of Step 6.
For every $t \in T$, let $Z_t$ be the critical locus of $l |_{F_t}$.  Note that
we have:
\beqn
Z_t = \{ \beta [i] (t) \}_{i \in I} \, .
\eeqn
For each $i \in I$ and each $t \in T$, consider the horizontal path:
\beqn
\gamma_h [i, t] = \gamma_h [(l \circ \beta [i]) (t)] : [0, 1] \to \bC,
\eeqn
as defined in Step 5.  Also, consider the positive eigenspace:
\beqn
T_h [i, t] \subset T_{\beta [i] (t)} F_t \, ,
\eeqn
defined by analogy with the positive eigenspace $T_h [i]$ of Step 6.
The real vector spaces $\{ T_h [i, t] \}_{t \in T}$ form a local system
over $T$.  Let $\cO_h [i, t]$ be the orientation of $T_h [i, t]$ obtained
from the orientation $\cO_h [i]$ of $T_h [i] = T_h [i, 0]$ of Step 6 via
the parallel transport of the local system $\{ T_h [i, t] \}_{t \in T}$.

{\em Step 14.}  The triple $(\beta [i] (t), \gamma_h [i, t], \cO_h [i, t])$
does not define a Picard-Lefschetz class in $M_t$ for every pair $(i, t)$.
This is because the path $\gamma_h [i, t]$ can collide with the set of
critical values $l (Z_t)$ for some $t > 0$.  However, by condition (B2) of
Step 10, such collisions can only occur for finitely many $t \in T$.  More
precisely, for each pair $i, j \in I$, let $\Theta_{i, j} \subset T$ be the set
of roots of the equation of condition (B2) of Step 10.  The set $\Theta_{i,j}$
is either empty or consists of one element.  Consider the union:
\beqn
\Theta = \bigcup_{i, j \in I} \Theta_{i, j} \, ;
\eeqn
it is a finite subset of $T$.  Define $T^\circ = T \setminus \Theta$.  Note
that, by condition (B2) of Step 10, we have $\{ 0, 1 \} \subset T^\circ$.
Now, for each $i \in I$ and each $t \in T^\circ$, we can use the triple
$(\beta [i] (t), \gamma_h [i, t], \cO_h [i, t])$ to define a Picard-Lefschets
class $u [i, t] \in M_t$ by analogy with the class $u [i]$ of Step 6.  Moreover,
just as in Step 8, for every $t \in T^\circ$, the elements $\{ u [i, t] \}_{i \in I}$
form a basis of $M_t$.

{\em Step 15.}  For each $i \in I$ and each pair $t_1, t_2 \in T^\circ$, such that 
$[t_1, t_2] \subset T^\circ$, we have:
\beqn
\nu_{t_1, t_2} (u [i, t_1]) = u [i, t_2].
\eeqn
Now, let $\theta \in \Theta$.  To describe what happens to the basis element
$u [i, t]$, as $t$ passes through $\theta$, we will need the following notation.
For each $i \in I$, let $J [i, \theta] \subset I$ be the set of all $j \in I$, such that:
\beqn
(\zeta \circ \beta [i]) (\theta) = (\zeta \circ \beta [j]) (\theta)
\;\; \text{and} \;\;
(\xi \circ \beta [i]) (\theta) < (\xi \circ \beta [j]) (\theta).
\eeqn
Let $t_1, t_2 \in T^\circ$ be a pair, such that: $t_1 < \theta < t_2$,
$[t_1, \theta) \subset T^\circ$, and $(\theta, t_2] \subset T^\circ$.
Then, by a standard Picard-Lefschetz theory argument, we have:
\beqn
\nu_{t_1, t_2} (u [i, t_1]) = u [i, t_2] +
\sum_{j \in J [i, \theta]} a_{i, j} \cdot u [j, t_2],
\eeqn
for some integers $\{ a_{i, j} \}_{j \in J [i, \theta]}$.

{\em Step 16.}  Recall the subspace $M^0 \subset M_0$, defined in Step 9.
For every $t \in T$, let $M^0_t = \nu_{0, t} (M^0) \subset M_t$.  Note that,
by the definition of $M^0$ and the last equation of Step 12, we have:
\beqn
M^0_1 = M^0_0 = M^0 \subset M_0 = M_1 \, .
\eeqn
The requisite containment $u [k] \in M^0$ will be implied by the following
claim.

\noindent
{\bf Claim:}  {\em 
Let $i \in I$ and $t \in T^\circ$.  Assume that:
\beqn
(\xi \circ \beta [i]) (t) > \xi (e_k).
\eeqn
Then we have $u [i, t] \in M^0_t$.
}

Before proving the Claim, we note that it readily implies that $u [k] \in M^0$.
Indeed, pick a $t \in T^\circ$, such that $[t, 1] \subset T^\circ$.  Consider
the element $u[0, t] \in M_t$.  By condition (B1) of Step 10 for $i = 0$,
we know that the composition $\xi \circ \beta [0] : [0, 1] \to \bR$ is monotone
decreasing.  Therefore, we have:
\beqn
(\xi \circ \beta [0]) (t) > (\xi \circ \beta [0]) (1) = \xi (e_k).
\eeqn
Thus, by the Claim, we have $u [0, t] \in M^0_t$.  Now, by the first equation
of Step 15, applied to $t_1 = t$ and $t_2 = 1$, we have:
\beqn
u [0, 1] = \nu_{t, 1} (u [0, t]) \in  \nu_{t, 1} (M^0_t) = M^0_1 = M^0_0
\subset M_0 \, .
\eeqn
It remains to note that, by construction, we have:
\beqn
u [k] = u [k, 0] = \pm u [0, 1] \in M_0 = M_1 \, .
\eeqn
Indeed, the Picard-Lefschetz classes $u [k, 0]$ and $u [0, 1]$ are defined
using the same critical point $e_k \in Z$ and the same path
$\gamma_h [k] : [0, 1] \to \bC$.  The only difference between these two
classes comes form the difference between the orientations $\cO_h [k, 0]$
and $\cO_h [0, 1]$ of the real vector space $T_h [0, k] = T_h [1, 0]$.
Therefore, these classes are the same up to sign.

Steps 17-22 below comprise our proof of the Claim.

{\em Step 17.}  We proceed to prove the Claim of Step 16, arguing by
contradiction.  Assume that the claim is false.  Consider the set:
\beqn
S = \{ (i, t) \in I \times T^\circ \; | \; u [i, t] \notin M^0_t \}.
\eeqn
Define:
\beqn
\xi_{max} = \sup_{(i, t) \in S} (\xi \circ \beta [i]) (t).
\eeqn
By assumption, we have $\xi_{max} > \xi (e_k)$.  

Let $C (T^\circ) = \pi_0 (T^\circ)$ be the set parametrizing the connected
components of $T^\circ$.  For each $\tau \in C (T^\circ)$, let $T^\circ [\tau]
\subset T^\circ$ be the corresponding connected component.  Thus, we 
have:
\beqn
T^\circ = \bigcup_{\tau \in C (T^\circ)} T^\circ [\tau].
\eeqn
The first equation of Step 15 implies that, for every pair $(i, \tau) \in
I \times C (T^\circ)$, we have either $(\{ i \} \times T^\circ [\tau]) \subset S$
or $(\{ i \} \times T^\circ [\tau]) \cap S = \emptyset$.  Define:
\beqn
C (S) = \{ (i, \tau) \in I \times C (T^\circ) \; | \;
(\{ i \} \times T^\circ [\tau]) \subset S \}.
\eeqn
Thus, we have:
\beqn
S = \bigcup_{ (i, \tau) \in C (S) } \{ i \} \times T^\circ [\tau].
\eeqn
For every $(i, \tau) \in C (S)$, let:
\beqn
\xi_{max} [i, \tau] = \sup_{t \in T^\circ [\tau]} (\xi \circ \beta [i]) (t).
\eeqn
Since $C (S)$ is a finite set, we have:
\beqn
\xi_{max} = \max_{(i, \tau) \in C (S)} \xi_{max} [i, \tau].
\eeqn
Therefore, there exists a pair $(i, \tau) \in C (S)$, such that
$\xi_{max} [i, \tau] = \xi_{max} > \xi (e_k)$.  Let us fix such a pair
$(i, \tau)$.  In subsequent steps, we analyze the possibilities for this 
pair.

{\em Step 18.}  We consider four cases for the pair $(i, \tau) \in C(S)$
fixed at the end of Step 17.  (Recall condition (B1) of Step 10.)

{\em Case 1:}  We have $(\xi \circ \beta [i])' (t) > 0$ for all $t \in T$,
and $1 \in T^\circ [\tau]$.

{\em Case 2:}  We have $(\xi \circ \beta [i])' (t) > 0$ for all $t \in T$,
and $1 \notin T^\circ [\tau]$.

{\em Case 3:}  We have $(\xi \circ \beta [i])' (t) < 0$ for all $t \in T$,
and $0 \in T^\circ [\tau]$.

{\em Case 4:}  We have $(\xi \circ \beta [i])' (t) < 0$ for all $t \in T$,
and $0 \notin T^\circ [\tau]$.

\noindent
In Steps 19-22 below, we take up these four cases one-by-one, reaching
a contradiction in each case.

{\em Step 19.}  Suppose Case 1 of Step 18 obtains.  Then we have:
\beqn
\xi_{max} = \xi_{max} [i, \tau] = (\xi \circ \beta [i]) (1).  
\eeqn
By the assumption of Step 17, we have $\xi_{max} > \xi (e_k)$.  Therefore,
by Step 4, we have $\beta [i] (1) = e_j$, for some $j < k$.  By the induction
hypothesis of Step 9, we have $u [j] \in M^0$.  But, by an argument as in 
Step 16, we have:
\beqn
u [j] = u [j, 0] = \pm u [i, 1] \in M_0 = M_1 \, .
\eeqn
Thus, we conclude that $u [i, 1] \in M^0 = M^0_1$ (see the first displayed
equation of Step 16), which contradicts the assumption that
$(i, \tau) \in C (S)$.

{\em Step 20.}  Suppose Case 2 of Step 18 obtains.  Let
$\theta = \sup (T^\circ [\tau])$.  Since $1 \notin T^\circ [\tau]$, we have
$\theta \in \Theta$.  We also have:
\beqn
\xi_{max} = \xi_{max} [i, \tau] = (\xi \circ \beta [i]) (\theta).
\eeqn
Recall the index set $J [i, \theta] \subset I$ of Step 15.  Pick a
$t_1 \in T^\circ [\tau]$.  Next, pick a $t_2 \in (\theta, 1]$, such that
$(\theta, t_2] \subset T^\circ$, and we have:
\beqn
(\xi \circ \beta [j]) (t_2) > \xi_{max} \, ,
\eeqn
for every $j \in J [i, \theta]$.  This is possible by the definition of $J [i, \theta]$.

Next, consider the final equation of Step 15 for the pair $t_1, t_2 \in T^\circ$.
We have:
\beqn
\nu_{t_1, t_2} (u [i, t_1]) = u [i, t_2] +
\sum_{j \in J [i, \theta]} a_{i, j} \cdot u [j, t_2].
\eeqn
By assumption, we have $u [i, t_1] \notin M^0_{t_1}$.  By the definition
of the subspaces $\{ M^0_t \subset M_t \}_{t \in T}$, we have
$\nu_{t_1, t_2} (M^0_{t_1}) = M^0_{t_2}$.  Therefore, we have
$\nu_{t_1, t_2} (u [i, t_1]) \notin M^0_{t_2}$.  On the other hand, every term
in the RHS of the above equation is contained in $M^0_{t_2}$, since by
construction, all of them correspond to critical points of $l |_{F_{t_2}}$ with
$\xi > \xi_{max}$.  Thus, we obtain a contradiction.

{\em Step 21.}  Suppose Case 3 of Step 18 obtains.  Then we have:
\beqn
\xi_{max} = \xi_{max} [i, \tau] = (\xi \circ \beta [i]) (0) = \xi (e_i).  
\eeqn
By the assumption of Step 17, we have $\xi_{max} > \xi (e_k)$.  Therefore,
by Step 4, we have $i < k$.  But, by the induction hypothesis of Step 9, we
have $u [i, 0] = u [i] \in M^0$, which contradicts the assumption that
$(i, \tau) \in C(S)$.

{\em Step 22.}  Suppose Case 4 of Step 18 obtains.  Our argument is
then parallel to Step 20.  Let $\theta = \inf (T^\circ [\tau])$.  Since
$0 \notin T^\circ [\tau]$, we have $\theta \in \Theta$.  We also have:
\beqn
\xi_{max} = \xi_{max} [i, \tau] = (\xi \circ \beta [i]) (\theta).
\eeqn
Recall the index set $J [i, \theta] \subset I$ of Step 15.  Pick a 
$t_1 \in [0, \theta)$, such that $[t_1, \theta) \subset T^\circ$.
Next, pick a $t_2 \in T^\circ [\tau]$, such that we have:
\beqn
(\xi \circ \beta [j]) (t_2) > \xi_{max} \, ,
\eeqn
for every $j \in J [i, \theta]$.  This is possible by the definition of $J [i, \theta]$.

Next, we can rewrite the final equation of Step 15 for the pair
$t_1, t_2 \in T^\circ$ as follows:
\beqn
u [i, t_2] = \nu_{t_1, t_2} (u [i, t_1]) -
\sum_{j \in J [i, \theta]} a_{i, j} \cdot u [j, t_2].
\eeqn
By assumption, we have $u [i, t_2] \notin M^0_{t_2}$.  On the other hand,
every term in the RHS of the above equation is contained in $M^0_{t_2}$.
Indeed, we have:
\beqn
(\xi \circ \beta [i]) (t_1) > (\xi \circ \beta [i]) (\theta) = \xi_{max} \, .
\eeqn
Therefore, we have $u [i, t_1] \in M_{t_1}$ and 
$\nu_{t_1, t_2} (u [i, t_1]) \in M_{t_2}$.  Finally, by the choice of $t_2$
above, we have $u [j, t_2] \in M^0_{t_2}$ for every $j \in J [i, \theta]$.
Thus, we obtain a contradiction, which completes our proof of
Lemma \ref{modified4.2}.

\subsection{Proof of Modified \cite[Lemma 4.3]{Gr3}}
\label{subsec-proof-4p3}

Our proof of Lemma \ref{modified4.3} is closely parallel to the proof of Lemma
\ref{modified4.2} described in the previous subsection.  In this subsection, we
only indicate the adaptations that are needed for Lemma \ref{modified4.3}.
Note that \cite[Lemma 4.3]{Gr3} is stated in terms of the local system $\cL$,
and the stalk $\cL_l$ of $\cL$ at $l$.  By definition, we have:
\beqn
\cL_l = H^{-d}_l (\cF \, P). 
\eeqn
In the argument below, we will mostly be referring to the stalk $\cL_l$ as the
RHS of the above equation, in order to keep the notation parallel with the proof 
of Lemma \ref{modified4.2} above.

Steps 1-7 of the proof of Lemma \ref{modified4.2} carry over to Lemma
\ref{modified4.3} without any changes.  In Step 8, we need to omit the last
sentence.  In Step 9, the first sentence needs to be modified as follows.

{\em Step 9.}  Let $M^0 \subset H^{-d}_l (\cF \, P)$ be the linear subspace
generated by the image of $u [0]$ under the holonomy action:
\beqn
h : \pi_1 ((V^*)^{rs}, l) \to \on{End} (H^{-d}_l (\cF \, P)).
\eeqn

\noindent
The rest of Step 9 is unchanged.  Steps 10 through 14 need to be modified
as follows.

{\em Step 10.}  Recall that we have $l \in c^*$, where $c^* \subset V^*$
is the Cartan subspace of \cite[Proposition 2.13]{Gr3}.  Consider the regular
part $(c^*)^{reg} = c^* \cap (V^*)^{rs}$.  For each $w \in W$, we will write
$w^* : c^* \to c^*$ for the adjoint operator of $w : c \to c$.  We have
the following claim.

\noindent
{\bf Claim:} {\em
There exists a smooth path:
\beqn
\beta : [0, 1] \to (c^*)^{reg},
\eeqn
with $\beta (0) = l$ and $\beta (1) = (w_k)^* \, l$, satisfying conditions
{\em (B1)} and {\em (B2)} below.  

{\em (B1)}  For each $i \in I$, the real part $\on{Re} (\beta' (t) (e_i))$
does not vanish at any point $t \in [0, 1]$.

{\em (B2)}  For each pair $i, j \in I$, there is at most one $t \in [0, 1]$, such that:
\beqn
\on{Im} (\beta (t) (e_i)) = \on{Im} (\beta (t) (e_j)).
\eeqn
Moreover, the root $t$ satisfies $t \in (0, 1)$ and:
\beqn
\on{Im} (\beta' (t) (e_i)) \neq \on{Im} (\beta' (t) (e_j)).
\eeqn
}

\noindent
{\bf Proof:} The path $\beta$ can be obtained as a $C^1$-small perturbation
of the straight line path $\beta_{st} : [0, 1] \to c^*$, given by:
\beqn
\beta_{st} : t \mapsto (1-t) \cdot l + t \cdot (w_k)^* \, l,
\eeqn
as in Step 10 of the proof of Lemma \ref{modified4.2}. 
\qed

{\em Step 11.}  Let $G \cdot l \subset (V^*)^{rs}$ be the $G$-orbit of $l$.  
Recall that $G$ is connected.  Therefore, so is $G \cdot l$.  Note that
$\beta (1) = (w_k)^* \, l \in G \cdot l$.  Pick a smooth path:
\beqn
\beta_1 : [0, 1] \to G \cdot l,
\eeqn
with $\beta_1 (0) = (w_k)^* \, l$ and $\beta_1 (1) = l$.  Define a loop:
\beqn
\bar \beta : [0, 1] \to (V^*)^{rs},
\eeqn
by composing the paths $\beta$ and $\beta_1$.  More precisely, let:
\beqn
\bar \beta (t) = \left\{ \begin{array}{ll}
\beta (2 \, t)         & \text{ for } t \in [0, 1/2] \\
\beta_1 (2 \, t -1) & \text{ for } t \in [1/2, 1].
\end{array}\right.
\eeqn
Note that  $\bar \beta$ represents an element of $\pi_1 ((V^*)^{rs}, l)$.  We will
show that $u [k] \subset M^0$ (see Step 9 above) by examining the holonomy
operator:
\beqn
h (\bar \beta) : H^{-d}_l (\cF \, P) \to H^{-d}_l (\cF \, P).
\eeqn

{\em Step 12.}  To analyze the holonomy operator $h (\bar \beta)$, we
need to consider a parametrized version of the basis $\{ u [i] \}_{i \in I}$ of
Step 8.  We begin by defining a parametrized version of the stalk cohomology
group $H^{-d}_l (\cF \, P)$.  Let us write $T = [0, 1]$.  For every $t \in T$, let:
\beqn
M_t = H^{-d}_{\beta (t)} (\cF \, P) = \cL_{\beta (t)} \, .
\eeqn
Note that the vector spaces $\{ M_t \}_{t \in T}$ form a local system over
the segment $T$.  For each pair $t_1, t_2 \in T$, let:
\beqn
\nu_{t_1, t_2} : M_{t_1} \to M_{t_2} \, ,
\eeqn
be the parallel transport of the local system $\{ M_t \}_{t \in T}$.  Also, let:
\beqn
\nu_{1, 2} : M_1 \to M_0 \, ,
\eeqn
be the holonomy of the local system $\cL$ along the path $\beta_1$.
Note that, with this notation, we have:
\beqn
\nu_{1, 2} \circ \nu_{0, 1} = h (\bar \beta) : M_0 \to M_0 \, .
\eeqn

{\em Step 13.}  
Next, we consider parametrized versions of the constructions of Step 6.
For every $t \in T$, let $Z_t$ be the critical locus of $\beta (t) |_F$.  Note that
the set $Z_t$ is independent of $t \in T$.  More precisely, for every $t \in T$,
we have:
\beqn
Z_t = Z = \{ e_i \}_{i \in I} \, .
\eeqn
For each $i \in I$ and each $t \in T$, consider the path:
\beqn
\gamma_h [i, t] = \gamma_h [\beta (t) (e_i)] : [0, 1] \to \bC,
\eeqn
as defined in Step 5.  Also, consider the positive eigenspace:
\beqn
T_h [i, t] \subset T_{e_i} F,
\eeqn
defined by analogy with the positive eigenspace $T_h [i]$ of Step 6, and
using $\beta (t) \in (c^*)^{reg}$ in place of $l$.  The real vector spaces
$\{ T_h [i, t] \}_{t \in T}$ form a local system over $T$.  Let $\cO_h [i, t]$ be
the orientation of $T_h [i, t]$ obtained from the orientation $\cO_h [i]$ of
$T_h [i] = T_h [i, 0]$ of Step 6 via the parallel transport of the local system
$\{ T_h [i, t] \}_{t \in T}$.

{\em Step 14.}  Just as in Step 14 of the proof of Lemma \ref{modified4.2},
we define a finite subset $\Theta \subset T$, and its complement
$T^\circ = T \setminus \Theta \supset \{ 0, 1 \}$.  Now, for each $i \in I$
and each $t \in T^\circ$, let $u [i, t] \in M_t$ be the Picard-Lefschets class
defined by the triple $(e_i, \gamma_h [i, t], \cO_h [i, t])$.  Just as in Step 8,
for every $t \in T^\circ$, the elements $\{ u [i, t] \}_{i \in I}$ form a basis of
$M_t$.

In Step 15, we only need to modify the definition of the index set
$J [i, \theta]$.  Namely, for each $i \in I$ and each $\theta \in \Theta$,
let $J [i, \theta] \subset I$ be the set of all $j \in I$, such that:
\beqn
\on{Im} (\beta (\theta) (e_i)) = \on{Im} (\beta (\theta) (e_j))
\;\; \text{and} \;\;
\on{Re} (\beta (\theta) (e_i)) < \on{Re} (\beta (\theta) (e_j)).
\eeqn
The rest of Step 15 is unchanged.  Step 16 needs to be modified as follows.

{\em Step 16.}  Recall the subspace $M^0 \subset M_0$, defined in Step 9.
For every $t \in T$, let $M^0_t = \nu_{0, t} (M^0) \subset M_t$.  Note that,
by the definition of $M^0$ and the last equation of Step 12, we have:
\beqn
M^0_0 = \nu_{1,2} (M^0_1) = M^0 \subset M_0.
\eeqn
The requisite containment $u [k] \in M^0$ will be implied by the following
claim.

\noindent
{\bf Claim:}  {\em 
Let $i \in I$ and $t \in T^\circ$.  Assume that:
\beqn
\on{Re} (\beta (t) (e_i)) > \xi (e_k).
\eeqn
Then we have $u [i, t] \in M^0_t$.
}

Before proving the Claim, we note that it readily implies that $u [k] \in M^0$.
Indeed, pick a $t \in T^\circ$, such that $[t, 1] \subset T^\circ$.  Consider
the element $u[0, t] \in M_t$.  By condition (B1) of Step 10 for $i = 0$,
we know that the real part:
\beqn
\on{Re} (\beta (t) (e_0)) : [0, 1] \to \bR,
\eeqn
is monotone decreasing.  Therefore, we have:
\beqn
\on{Re} (\beta (t) (e_0)) > \on{Re} (\beta (1) (e_0)) = \xi (e_k),
\eeqn
where the last equality follows from the condition $\beta (1) = (w_k)^* \, l$
of Step 10.  Thus, by the Claim, we have $u [0, t] \in M^0_t$.  Now, by the
first equation of Step 15, applied to $t_1 = t$ and $t_2 = 1$, we have:
\beqn
u [0, 1] = \nu_{t, 1} (u [0, t]) \in \nu_{t, 1} (M^0_t) = M^0_1 \, .
\eeqn
It remains to note that, by construction, and by the first equation of this
step, we have:
\beqn
u [k] = u [k, 0] = \pm \nu_{1,2} (u [0, 1]) \in M^0_0 \, .
\eeqn
More precisely, for the second equality above, observe that the set of
critical values of $\beta_1 (t) |_F$ is independent of $t$ for $t \in [0, 1]$.
The Picard-Lefschetz classes $u [k, 0]$ and $u [0, 1]$ correspond to
critical points of $\beta_1 (0) |_F$ and $\beta_1 (1) |_F$, respectively,
with the same critical value $l (e_k)$.  Furthermore, the classes $u [k, 0]$
and $u [0, 1]$ are defined using the same path:
\beqn
\gamma_h [k, 0] = \gamma_h [0, 1] = \gamma_h [l (e_k)] : [0, 1] \to \bC.
\eeqn
It follows that the classes $u [k, 0]$ and $\nu_{1,2} (u [0, 1])$ correspond
to the same critical point $e_k$ of $l |_F$, and the same path
$\gamma_h [l (e_k)] : [0, 1] \to \bC$.  Therefore, these classes are
the same up to sign.

Our proof of the Claim of Step 16 is a straightforward adaptation of Steps
17-22 of the proof of Lemma \ref{modified4.2}.  We omit the remaining
details.  This completes our proof of Lemma \ref{modified4.3}.

\subsection{Further Comments on the Proof of \cite[Theorem 3.1]{Gr3}:
the Stable Case}
\label{subsec-stable-case}

The proof of \cite[Theorem 3.1]{Gr3} in the stable case is readily fixed by
replacing \cite[Lemmas 4.2 $\&$ 4.3]{Gr3} by Lemmas \ref{modified4.2}
and \ref{modified4.3} above, and by further making the following four
changes.  Note that two out of these changes (the first and the third) are
unrelated to the main problem described in Section \ref{subsec-problem},
as we take this opportunity to correct some unrelated inaccuracies.

First, in stating \cite[Lemma 2.9]{Gr3}, we neglected to require that
$X \subset \hat X$ be a stratum of the stratification $\hat \cX$ of 
$\hat X$ provided by that lemma.  The construction of
[Gr1, Section 3.4] provides a stratification $\hat \cX$ satisfying this
requirement.  Moreover, the lemma with this requirement is trivially
implied by the lemma as stated.  Thus, in applying
\cite[Lemma 2.9]{Gr3} in the second paragraph of 
\cite[Section 4.2]{Gr3}, we must assume that $F$ is a stratum of
the stratification of $\hat F$.  This justifies the conclusion that 
$Z \subset c \subset V$, where $Z$ is the stratified critical locus of
the restriction $\hat l |_{\hat F}$.

Second, at \cite[p. 424]{Gr3}, we must choose $e_0$ and $\gamma_0$ so 
that conditions (A1) and (A2) are satisfied.

Third, in the second paragraph of the proof of claim (iii) of the theorem
at \cite[p. 425]{Gr3}, we must use the homomorphism
$\eta^{}_\lambda : B_W (\lambda) \to W^\sharp (Z)$
of Section \ref{subsec-braid-group} above, to replace the sentence:

\noindent
``Each cluster consists of $n$ points surrounding a point in $W \cdot v_1$,
and the action of $\sigma \in W$ cyclically permutes the points in each 
cluster.''

\noindent
by:

\noindent
``Each cluster consists of $n$ points surrounding a point in $W \cdot v_1$,
and the action of $\eta^{}_\lambda (\sigma) \in W^\sharp (Z)$ cyclically
permutes the points in each cluster.''

And fourth, in last paragraph of the proof of claim (iii) of the theorem, we
must modify the construction of the Picard-Lefschetz class $u_0$, so that 
Lemma \ref{modified4.3} would apply.  Namely, we must choose a point
$e_0 \in Z$, satisfying condition (A1).  Next, we must choose a smooth path
$\gamma_0 : [0, 1] \to \bC$, satisfying conditions (i)-(iv) of \cite[p. 424]{Gr3}
for $e = e_0$, as well as condition (A2).  Finally, we must pick an orinetation
$\cO_0$ of the real subspace $T_{e_0} [\gamma_0] \subset T_{e_0} F$, 
and define $u_0 = u (e_0, \gamma_0, \cO_0)$.  Figure 1 at \cite[p. 426]{Gr3}
must be modified accordingly, to show the path $\gamma_0$ originating from
$l (e_0)$, one of the right-most points of the image $l (Z)$.  The rest of the 
argument in the stable case remains unchanged.

\subsection{Further Comments on the Proof of \cite[Theorem 3.1]{Gr3}:
the Nonstable Case}
\label{subsec-nonstable-case}

As indicated in \cite[Section 4.3]{Gr3}, the proof of \cite[Theorem 3.1]{Gr3}
in the case when $G | V$ is not stable is analogous to the proof in the stable
case, with Morse critical points of the restriction $l |_F$ replaced by Morse-Bott
critical manifolds (see \cite[Corollary 2.16 (ii)]{Gr3}.  In this subsection, we
indicate the changes to the proof that have to be made in the nonstable case.
Compared to \cite[Section 4.3]{Gr3}, we reflect the changes described in
Section \ref{subsec-stable-case}, and supply some more detail.

As in the stable case, we pick a basepoint $l \in (V^*)^{rs}$, lying in the Cartan
subspace $c^*$ of \cite[Proposition 2.13]{Gr3}, consider the compactification
$\hat F$ of $F$ relative to $l$, fix a stratification of $\hat F$ as in
[Gr3, Lemma 2.9], such that $F \subset \hat F$ is a stratum, and write $Z$
for the stratified critical locus of the restriction $\hat l |_{\hat F}$.  Let
$Z_1 = Z \cap c$, it is a single $W$-orbit in $c$.  Note that, by 
[Gr3, Corollary 2.16 (ii)], we have:
\beq\label{eqn-Z-cap-F}
l (Z \cap F) = l (Z_1),
\eeq
but the image $\hat l (Z)$ is potentially larger than $l (Z_1)$.

The construction of Picard-Lefschetz classes described at [Gr3, p. 424] can
be adapted to this new setting as follows.  Let $e \in Z_1$ and let $\gamma :
[0, 1] \to \bC$ be a smooth path satisfying conditions (i)-(iv) of
\cite[p. 424]{Gr3}.  Unlike in the stable case, the Hessian $\cH_e :
T_e F \to \bC$ of $l |_F$ at $e$ is now degenerate.  However, we may still
consider the positive eigenspace $T_e [\gamma] \subset T_e F$ of the real
quadratic form $\on{Re} (\cH_e / \gamma' (0))$.  By [Gr3, Corollary 2.16 (ii)],
we have $\dim_\bR \, T_e [\gamma] = d_0 - r$.  Fix an orientation $\cO$ of
$T_e [\gamma]$.  The triple $(e, \gamma, \cO)$ defines a Picard-Lefschetz 
class:
\beqn
u = u (e, \gamma, \cO) \in H_{d_0 - r} (F, \{ \xi (y) \geq \xi_0 \}; \bC),
\eeqn
where $\xi = \on{Re} (l)$ and $\xi_0$ is large, exactly as in the stable
case.  Moreover, as before, we can use [Gr3, Lemma 2.8] and duality with
supports to regard $u$ as an element of $H^{-d_0}_l (\cF \, P)$.

In order to analyze such Picard-Lefschetz classes, we introduce the following
notation.  For each $x \in \hat l (Z)$, let $Z_x = Z \cap \hat l ^{-1} (x)$ and let:
\beqn
\hat l [x] = \hat l |_{\hat F} - x : \hat F \to \bC.
\eeqn
Consider the vanishing cycles sheaf:
\beqn
Q_x = \phi_{\hat l [x]} (j_! \, \bC_F \, [d - r]);
\eeqn
it is a a perverse sheaf on $\hat F$, satisfying:
\beq\label{eqn-support-of-Q}
\on{supp} (Q_x) \subset Z_x \, .
\eeq
By [Gr3, Lemma 2.8], we have:
\beq\label{eqn-stalk-as-direct-sum}
H^{-d_0}_l (\cF \, P) \cong \bigoplus_{x \in \hat l (Z)}
\bH^{d - d_0} (\hat F; Q_x),
\eeq
where the isomorphism depends on a choice of a system of cuts, as usual.

By the inequalities of [Gr3, Lemmas 2.9 $\&$ 4.1] and equations
\eqref{eqn-Z-cap-F}-\eqref{eqn-support-of-Q} above, we have:
\beq\label{eqn-invisible}
\bH^{d - d_0} (\hat F; Q_x) = 0 \;\; \text{for every} \;\;
x \in \hat l (Z) \setminus l (Z_1).
\eeq
In other words, the critical values $x \in \hat l (Z) \setminus l (Z_1)$ are
``invisible'' from the point of view of the direct sum decomposition
\eqref{eqn-stalk-as-direct-sum}.  This observation enables us to extend
the notation for Picard-Lefschetz classes as follows.  Let $e \in Z_1$ and
let $\gamma : [0, 1] \to \bC$ be a smooth path satisfying conditions (i),
(iii), and (iv) of \cite[p. 424]{Gr3}, plus the following condition:

(ii$'$) $\gamma (t) \notin l (Z_1)$, for $t > 0$.

\noindent
Define the real subspace $T_e [\gamma] \subset T_e F$ as usual, and
fix an orientation $\cO$ of $T_e [\gamma]$.  Then we can define:
\beqn
u (e, \gamma, \cO) = u (e, \gamma_1, \cO) \in H^{-d_0}_l (\cF \, P),
\eeqn
where $\gamma_1 : [0, 1] \to \bC$ is a smooth $C^1$-perturbation of the
path $\gamma$, satisfying conditions (i)-(iv) of \cite[p. 424]{Gr3}, and
such that $\gamma_1' (0) = \gamma' (0)$.  Equation \eqref{eqn-invisible}
ensures that the class $u (e, \gamma, \cO)$ is well defined.  This
enhances the analogy with the stable case, as we can basically ignore
the critical values $x \in \hat l (Z) \setminus l (Z_1)$ for the purposes of
discussing the paths used to define Picard-Lefschetz classes.

The main distinction from the stable case is that there is no a priori
geometric reason for the class $u (e, \gamma, \cO)$ to be non-zero.
This is because the Morse-Bott manifold containing $e$ is non-compact.
Thus, we can not draw an immediate conclusion from equation
\eqref{eqn-stalk-as-direct-sum} regarding the dimension
$\dim \, H^{-d_0}_l (\cF \, P)$.  To analyze this dimension, let:
\beqn
(c^*)^{rs, \circ} [\lambda] = \{ l' \in (c^*)^{rs} \;\, | \;\, |l' (Z_1)| = |Z_1| = |W| \}.
\eeqn
Note that $(c^*)^{rs, \circ} [\lambda]$ is Zariski open in $(c^*)^{rs}$.
We incorporate the regular vlalue $\lambda \in Q^{reg}$ in the notation
to indicate the dependence on the fiber $F = f^{-1} (\lambda)$.

Assume that $l \in (c^*)^{rs, \circ} [\lambda]$.  From the inequalities of
[Gr3, Lemmas 2.9 $\&$ 4.1], [Gr3, Corollary 2.16 (ii)], and equation
\eqref{eqn-support-of-Q} above, we can conclude that:
\beqn
\dim \, \bH^{d - d_0} (\hat F; Q_x) \in \{ 0, 1 \} \;\; \text{for every} \;\;
x \in l (Z_1).
\eeqn
Moreover, by construction, we have:
\beq\label{eqn-vanishing-equivalence}
u (e, \gamma, \cO) = 0 \iff \bH^{d - d_0} (\hat F; Q_x) = 0,
\eeq
where $x = l (e)$.  We can use the equivalence
\eqref{eqn-vanishing-equivalence} and the monodromy action of $B_W$ 
on $H^{-d_0}_l (\cF \, P)$ to conclude that if the hypercohomology group
$\bH^{d - d_0} (\hat F; Q_x)$ vanishes for some $x \in l (Z_1)$, then it
vanishes for all $x \in l (Z_1)$.  However, in the latter case, by
\eqref{eqn-stalk-as-direct-sum}, it would follow that
$H^{-d_0}_l (\cF \, P) = 0$.  Since the same conclusion would hold for
every $l' \in (c^*)^{rs, \circ} [\lambda]$, by [Gr3, Theorem 2.6] and
[Gr3, Corollary 2.16 (i)], it would then follow that $P = 0$.  This
contradiction shows that:
\beq\label{eqn-dim-one-term}
\dim \, \bH^{d - d_0} (\hat F; Q_x) = 1 \;\; \text{for every} \;\; x \in l (Z_1),
\eeq
and we have:
\beq\label{eqn-stalk-dimension}
\dim \, H^{-d_0}_l (\cF \, P) = |W|.
\eeq

So far, we have only established equation \eqref{eqn-stalk-dimension}
for $l \in (c^*)^{rs, \circ} [\lambda]$.  However, the LHS of this equation
is independent of the choice of the regular value $\lambda \in Q^{reg}$,
and the sets $(c^*)^{rs, \circ} [\lambda]$, for different choices of 
$\lambda$, cover the entirety of $c^{rs}$.  Therefore, equation
\eqref{eqn-stalk-dimension} holds for all $l \in (c^*)^{rs}$, and by 
$G$-equivariance, for all $l \in (V^*)^{rs}$.  Putting this together with
[Gr3, Theorem 2.6] and [Gr3, Corollary 2.16 (i)], we can conclude that:
\beqn
\cF \, P |_{(V^*)^{rs}} \cong \cL [d_0],
\eeqn
where $\cL$ is a rank $|W|$ local system on $(V^*)^{rs}$.

Note that equations
\eqref{eqn-vanishing-equivalence}-\eqref{eqn-dim-one-term}
imply that $u (e, \gamma, \cO) \neq 0$ whenever
$l \in (c^*)^{rs, \circ} [\lambda]$.  The same conclusion can be readily
extended to all $l \in (c^*)^{rs}$, for example, by considering the parallel
transport of the local system $\cL$ from a point $l \in (c^*)^{rs} \setminus
(c^*)^{rs, \circ} [\lambda]$ to a nearby point $l_1 \in (c^*)^{rs, \circ}
[\lambda]$.  Thus, we see that the Picard-Lefschetz class
$u (e, \gamma, \cO) \in H^{-d_0}_l (\cF \, P)$ is non-zero for every triple
$(e, \gamma, \cO)$, as in the stable case.

The rest of the proof of [Gr3, Theorem 3.1] in the nonstable case
proceeds by working with such Picard-Lefschetz classes in complete
analogy with the stable case.  We omit the rest of the details.

\section{Comments on \cite[Sections 5-6]{Gr3}}
\label{sec-sec-5-6}

In this section, we correct several minor errors from
\cite[Sections 5-6]{Gr3}.

First, in the second paragraph of \cite[Section 5]{Gr3}, the subgroup
$G_\sigma \subset G$ should be defined as the connected subgroup
corresponding to the subalgebra $\Lg_\sigma \subset \Lg$, not as the
adjoint form of $\Lg_\sigma$.

Second, in the last paragraph of the proof of \cite[Theorem 5.2]{Gr3},
the first sentence should be modified as follows:

\noindent
``Let $g_\sigma : Q_\sigma \to Q$ be the map
$g_\sigma : f_\sigma (v) \mapsto f (v + v_1)$, $v \in V_\sigma$.''

\noindent
And in the next line, the words ``the point $f_\sigma (v_1)$'' should
be replaced by ``the origin.''

Third, in the first paragraph of \cite[Section 6]{Gr3}, the sentence:

\noindent
``Then $\Lg^+$ is a Lie algebra, and the adjoint form $G^+$ of
$\Lg^+$ acts on the symmetric space $\Lg^-$ by conjugation.''

\noindent
should be replaced by:

\noindent
``Then $\Lg^+$ is a Lie algebra, and if $G$ is the adjoint for of $\Lg$,
and $G^+ \subset G$ is the connected subgroup corresponding
to $\Lg^+ \subset \Lg$, then $G^+$ acts on the symmetric space
$\Lg^-$ by conjugation.''

And fourth, in the proof of \cite[Theorem 6.1 (i)]{Gr3}, the sentence:

\noindent
``It follows from the root space decomposition for $\Lg_\bR$ that, 
in terms of Proposition 5.1, we have $G_\sigma = SO (s (i) + 1)$,
$V_\sigma = \bC^{s (i) + r}$, and $G_\sigma$ acts by the standard
representation on the first $s (i) + 1$ coordinates.''

\noindent
should be replaced by:

\noindent
``It follows from the root space decomposition for $\Lg_\bR$ that, 
in the notation of Section 5, there are isomorphisms
$V_\sigma \cong \bC^{s (i) + 1} \oplus c_\sigma$ and
$Q_\sigma \cong \bC \oplus c_\sigma$, such that
$f_\sigma = (\bar f_\sigma, \on{Id}_{c_\sigma})$, where
$\bar f_\sigma : \bC^{s (i) + 1} \to \bC$ is a non-degenerate
homogenous quadric.''

We conclude this section with a remark regarding the proof of
\cite[Theorem 6.1 (ii)]{Gr3}.

\begin{remark}
\label{rmk-gvx}
The proof of \cite[Theorem 6.1 (ii)]{Gr3} is essentially left as an exercise
to the reader.  Arguably, more should have been said.  It is indeed easy
to check that the answer provided for the homomorphism $\rho$ is
accurate up to sign.  One way to pin down the signs is to trace the effect
of both the monodromy action $\mu$ and the holonomy of the local system
$\cL$ on the top homology group $H_{d-r} (F; \bC) \cong \bC$, regarded
as a subgroup of the stalk $\cL_l$ (see the first displayed equation at
\cite[p. 424]{Gr3}).

The recent paper \cite{GVX} can be regarded as a sequel to
\cite[Section 6]{Gr3}, giving a generalization of \cite[Theorem 6.1 (ii)]{Gr3}
to nearby cycles with coefficients in certain rank one local systems on
the general fiber $F$ of $f$.  In particular, \cite[Section 6]{GVX} treats
the case of constant coefficients in detail, and \cite[Proposition 6.14]{GVX}
provides the analysis of the top homology group $H_{d-r} (F; \bC)$
mentioned above.  The relationship between the material in
\cite[Section 6]{Gr3} and \cite[Section 6]{GVX} is explained in the
introduction to the latter.
\end{remark}

\end{document}